# MISSION ANALYSIS FOR THE ION BEAM DEFLECTION OF FICTITIOUS ASTEROID 2015PDC

**Claudio Bombardelli**[(1)], **Davide Amato**[(2)], **Juan Luis Cano**[(3)]

[(1)(2)]*Space Dynamics Group, Technical University of Madrid, Plaza Cardenal Cisneros 3, 28040 Madrid, +34 913363939, claudio.bombardelli@upm.es, d.amato@upm.es*

[(3)]*Elecnor Deimos, Ronda de Poniente, 19, 28760 Tres Cantos, Madrid, juan-luis.cano@deimos-space.com*

***Keywords:*** *2015 PDC, ion beam shepherd, deflection, rendezvous, orbit uncertainty*

## ABSTRACT

A realistic mission scenario for the deflection of fictitious asteroid 2015PDC is investigated that makes use of the ion beam shepherd concept as primary deflection technique. The article deals with the design of a low thrust rendezvous trajectory to the asteroid, the estimation of the propagated covariance ellipsoid and the outcome of a slow-push deflection starting from three worst case scenarios (impacts in New Delhi, Dhaka and Teheran). Displacing the impact point towards very low populated areas, as opposed to full deflection, is found to be the simplest and most effective mitigation approach. Mission design, technical and political aspects are discussed.

## 1- INTRODUCTION

A hypothetical asteroid impact scenario was presented at the 2015 planetary defense conference (PDC) in Frascati, Rome with the aim of stimulating the discussion on several aspects of asteroid threat mitigation. The fictitious 150-500 m diameter asteroid, discovered on April 13 2015 and named 2015 PDC, was found to have several potential impacts with the Earth, the earliest and most likely on September 3, 2022. The impact probability, estimated in mid-June 2015, would reach 1% and would continue to rise with the rest of the scenario to be played out at the conference.

Starting from the asteroid ephemerides provided by the conference organizers and propagating forward until impact with the Earth one finds a nominal impact point in the South China Sea roughly 550 km off the Vietnam coast at around 03:52:10 UT, with about 16 km/s impact velocity at roughly 56 degrees from the surface. However, by considering a major line of variation error in the initial asteroid orbital determination one would obtain a *path of risk* stretching from eastern Turkey until the middle of the Pacific Ocean (almost 2000 km off the coast of Mexico) and passing through heavily populated areas such as Northern India, and major cities like New Delhi, Tehran and Dhaka.

Owing to the asteroid orbit geometry with respect to the Earth there appear to be no significant opportunities to increase the accuracy of the asteroid orbit by future ground-based observations.

Given the potentially devastating outcome of this scenario the first action to be taken is, arguably, to send a spacecraft to rendezvous with the asteroid and drastically reduce the uncertainty of its orbit by "collaborative" orbit determination [1]. We will show that, starting from a reasonable expected level of orbit determination accuracy, the size of the 1-sigma impact ellipse for 2015PDC would reduce to about 100 km at the predicted impact date, which can provide a quite definite answer to whether (and, if yes, where) an impact would take place. In addition, the rendezvous will provide an estimate of the asteroid mass to within a few percent error (as done for Itokawa).

In the case an impact with the Earth is predicted to occur the first thing to do will be to carefully assess its consequences. An impact in the open sea, for instance, might not be critical enough to warrant deflection as asteroids smaller than 500 m lack the capability of generating tsunami-like impact waves [2]. On the other hand, impacts to within about 100 km from coastlines or densely populated area would have catastrophic consequences [2], and the expected infrastructure damage including massive evacuation will most definitely dwarf the cost of a large space deflection mission.

This paper deals with three worst case scenarios in which the asteroid is predicted to impact in the vicinity of New Delhi (India), Dhaka (Bangladesh) and Tehran (Iran), and the design of a possible contactless ion beam deflection mission, a concept recently proposed by one of these authors [3,4]. The concept has the key advantage of exploiting high specific impulse ionic thrusters as an efficient mean to both transfer the spacecraft to a rendezvous trajectory with the asteroid and (if required) to produce an accurate deflection that can be tailored according to specific requirements.

The structure of the article is as follows. First, we analyze the asteroid impact path of risk and provide the reference orbital elements for the three impact scenarios previously described. Second, we design an interplanetary rendezvous trajectory to the asteroid starting from a preliminary propulsion system design and estimated power and mass budget. From the obtained rendezvous date, starting from an Itokawa-like covariance matrix, and by linear propagation up to the impact date an estimated impact covariance ellipse on the Earth surface is computed. Next, a deflection action using a continuous tangential thrust of variable magnitude according to power availability is implemented and the shift in the impact spot on the Earth surface is computed as a function of the thrust duration. Finally, an estimation of the impact damage variation as a function of the impact spot displacement is assessed by integrating a grey-scale nighttime picture of the Earth and georeferenced population density data over an expected circular damage zone of 100 km radius. Political issues related to possible crossing of the impact spot through neighboring countries and associated risks are discussed also in relation with the deflection strategy to be eventually chosen.

## 2- IMPACT SCENARIO

The reference equinoctial orbital elements of fictitious asteroid 2015PDC are provided in Table 1 and referring to the epoch of April 13 2015 at 0:00:00 UT.
The asteroid orbit is propagated up to the instant of intersection with the Earth ellipsoid using a simplified solar system model including gravitational perturbations from the solar system planets, the Moon (through JPL DE405 ephemerides) and the three largest asteroids[1] (Ceres, Vesta and Pallas) as well as first order relativistic corrections.

**Table 1**. Nominal ephemerides of fictitious asteroid 2015 PDC

| epoch (MJD) | 57125 |
|---|---|
| a (AU) | 1.775998173759480 |
| P1() | -0.448551534990503 |
| P2() | 0.198239860639469 |
| Q1() | -0.015660086557340 |
| Q2() | 0.043990645962994 |
| ML(deg) | 264.0060035482113 |

While the propagation of the nominal asteroid ephemerides would lead to an impact in the South China Sea (113.9 E 13.8 N) a preliminary path of risk estimate can be obtained by propagation starting from a set of slightly modified reference epochs. This would approximate a line of variation (LOV) sampling of the asteroid close to the direction of greatest orbit determination uncertainty [5].

The path of risk obtained following the above procedure is plotted in Figure 1 and stretches from Easter Turkey (37.9 N, 37.3 E) up to the South Pacific (6.7 N, -118.9 E) passing through heavily populated Earth regions including Northern India.

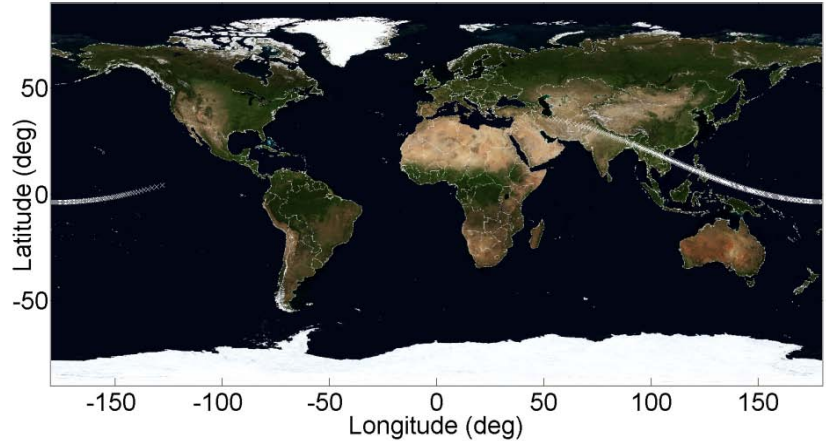

**Figure 1**. Path of risk of fictitious asteroid 2015PDC

---

[1]Note that, due to the eccentric orbit of 2015PDC, the perturbation from the largest asteroid belt objects can have an important effect. For high fidelity propagation it is recommended to employ as many as the largest objects as possible.

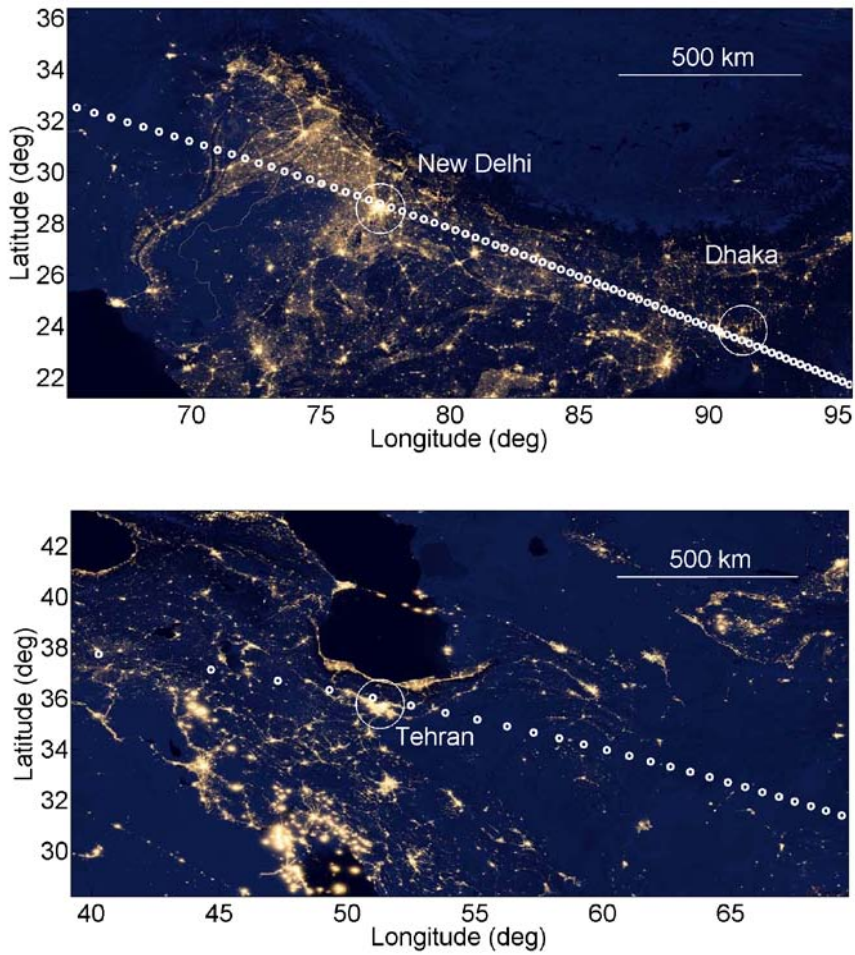

**Figure2**. Path of risk of fictitious asteroid 2015 PDC near New Delhi, Dhaka and Teheran

In particular, the three heavily populated capitals New Delhi, Dhaka and Tehran happen to lay very close to the path. The orbit initial conditions leading to an impact close to these cities are reported in Table 2 and were used as reference worst case scenarios for the deflection mission design conducted in this paper.

**Table 2**. Modified reference epoch of Table1 leading to worst case impacts

| epoch (MJD)    | Impact point coordinates | Nearest city |
|----------------|--------------------------|--------------|
| 57124.9984027  | 28.8 N, 77.2 E           | New Delhi    |
| 57124.9988310  | 23.8 N, 90.3 E           | Dhaka        |
| 57124.9979780  | 35.9 N, 51.4 E           | Tehran       |

## 2- RENDEZVOUS MISSION DESIGN

Given the quite high collision probability estimated in June 2015 and the lack of significant observation opportunities until the impact date it will be compulsory to design and launch a space mission to rendezvous with 2015PDC as soon as possible, track its trajectory and determine its size.
The mission will need to estimate the asteroid position and velocity to the greatest possible accuracy in order to confirm or rule out an impact. If an impact is indeed confirmed it will be important to know where and discuss possible countermeasures.

A first and very important tradeoff in the design of the rendezvous mission is whether the spacecraft should have a deflection capability, which would come at a price of higher complexity and weight.
At first sight, it does not look reasonable to embark a nuclear deflection device on the spacecraft. It would greatly complicate the mission and slow down its development. In addition, a nuclear deflection mission could be implemented at a later stage as the huge amount of transmitted momentum would probably allow, unlike other methods, for an effective deflection even if conducted just a few months before the expected impact.
On the other hand a slow-push deflecting capability based on contactless ion beam shepherding (IBS) would fit in a rendezvous mission in a much more comfortable way: Ion thrusters can be efficiently exploited for both the interplanetary transfer trajectory and the contactless deflection with no major expected increase in technological complexity. Finally, the cost of the added propellant mass required for deflection can be largely justified by the criticality of a possible impact event. Ultimately, for the mission to be analyzed in this article, the spacecraft will have an IBS deflection system on board with enough propellant to accomplish the deflection of a medium size asteroid, as it will be shown later.
Since the asteroid needs to be accurately tracked as soon as possible and enough time has to be given for an effective deflection the interplanetary mission should be designed with the fastest possible trajectory. With a high-specific-impulse low-thrust interplanetary trajectory in mind this means power requirements will be quite demanding. As a preliminary design we propose a Dawn-type 10kW(at 1AU) power subsystem with an estimated 150 kg mass but with two redundant sets of ionic thrusters mounted at a relative orientation of 180 degrees to provide ion beam shepherding capability and with 3500 s specific impulse, 70% electrical efficiency, and a 200mN + 200mN maximum thrust capability at 1AU. Following the Dawn spacecraft design each set of thrusters could contain 3 redundant units. We assume that a peak 400 mN thrust (1AU) can be reached during the interplanetary transfer by employing all available thruster-dedicated power to feed one set of thrusters. A preliminary estimate for the total spacecraft mass at interplanetary orbit insertion is of 1200 kg including 500 kg of Xenon (200 kg for the interplanetary trajectory, 300 kg for the deflection).

A low thrust trajectory optimization for the interplanetary phase has been performed assuming a launch between October 2015 and June 2017 and a Soyuz injection with a C3 up to 10 km$^2$/s$^2$ for our mass.

The optimization process assumes that the available power depends on the distance to the Sun with an inverse exponent of 1.7, instead of purely 2 (this is so to take into account that most solar cells gain performance at lower temperatures, thus at larger distances).

The most favorable trajectory found consists of a launch on May 28 2017 and arrival at the asteroid on Sep 30 2019 with a total of 200 kg of fuel spent. The transfer trajectory consists of a thrust-coast-thrust structure of 2.34 years of total duration. First thrust arc has a duration of 223 days, followed by a coasting phase of 412 days and a final thrust arc of 219 days. Figure 3 depicts the Ecliptic projection of the trajectory where the thrust arcs are represented with a thicker line. The thrust in the first arc is almost in the direction of the S/C velocity, whereas the second is closer to the opposite of the velocity.

Departure excess velocity is at the maximum possible for the available mass, i.e. 3.16 km/s and the departure declination is -9.1 deg, which would allow using the launcher from Kourou.

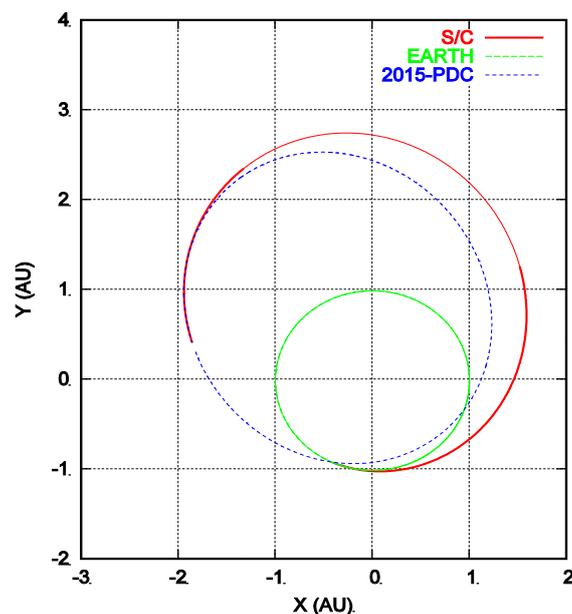

**Figure3**. Ecliptic projection of the low-thrust transfer trajectory from Earth to 2015 PDC launched in 2017

Figure 4 provides the S/C distances to the main bodies. Maximum distance to the Sun is below 2.8 AU, whereas maximum Earth distance is 3.6 AU. Distance to the Sun at arrival is 1.9 AU and to the Earth 2.8 AU. Figure 5 gives the S/C angles to the Sun and the Earth. Superior conjunctions occur in June 2018 and in September 2019, just before arrival to the asteroid. Under this condition it might be necessary to delay the arrival to the asteroid to have a good radiometric link to the S/C from Earth and thus a proper orbit determination solution. From the moment of arrival the visibility conditions from Earth shall be good (increasing sun elongation), with decreasing distances to Earth.

Figure 6 provides the evolution of the S/C mass, which reduces from 1200 kg to 1059 kg after the first thrust arc and subsequently from that value to 1006 kg at the end of the second thrust arc. Figure 7 provides the evolution of the thrust magnitude, which is dictated by the distance to the Sun with the previously commented dependency. Finally, Figure 8 provides the evolution of the thrust angle to directions of interest as the Sun-S/C line, the Earth-S/C line and the orbit normal.

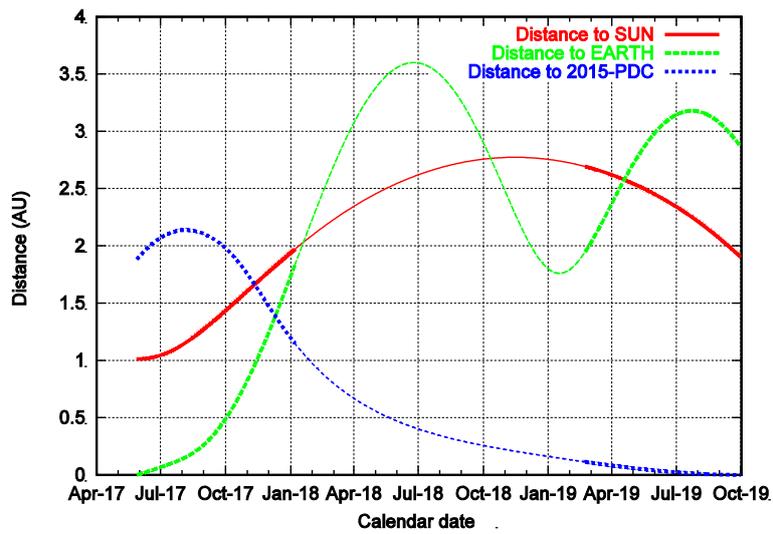
**Figure4**. Timely evolution of the S/C distances to the Sun, the Earth and 2015 PDC

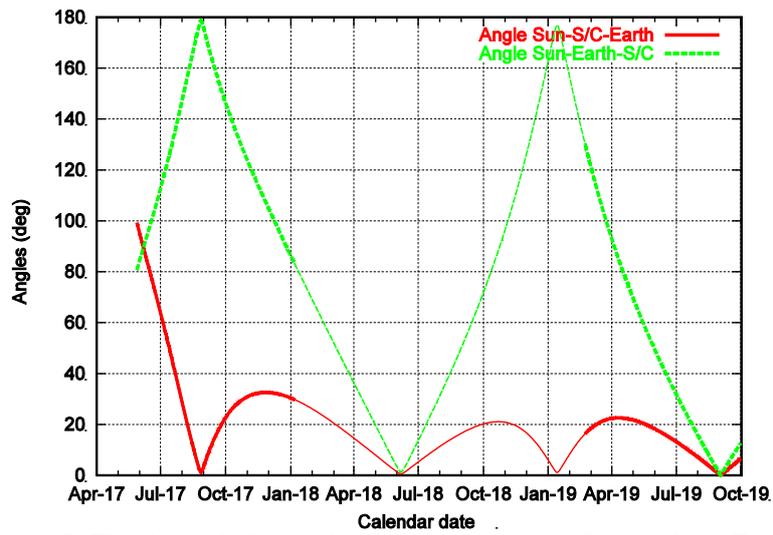
**Figure5**. Timely evolution of the S/C angles to the Sun and the Earth

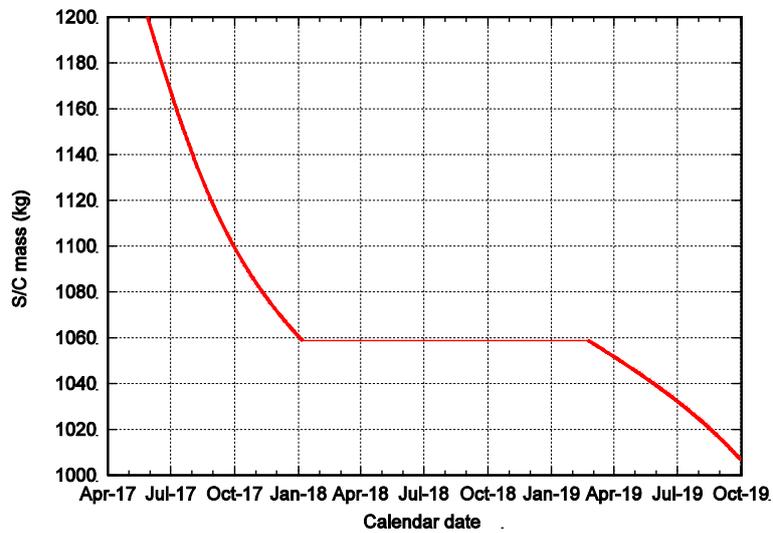
**Figure6**. Timely evolution of the S/C mass

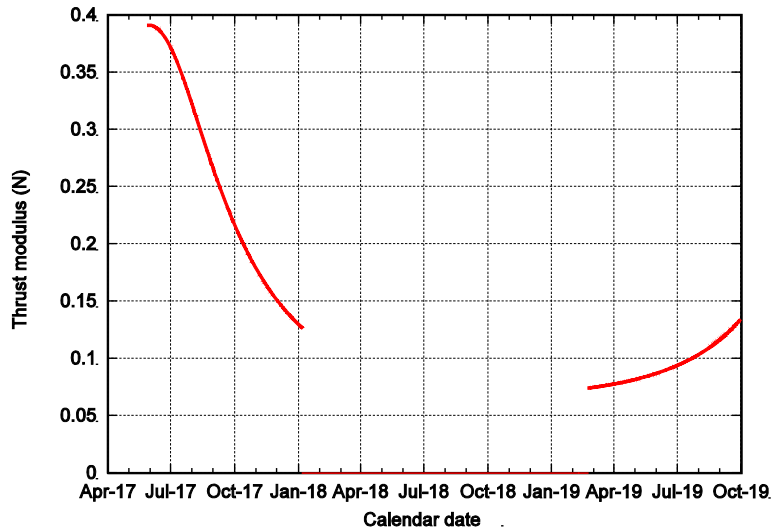
**Figure7**. Time evolution of the thrust magnitude

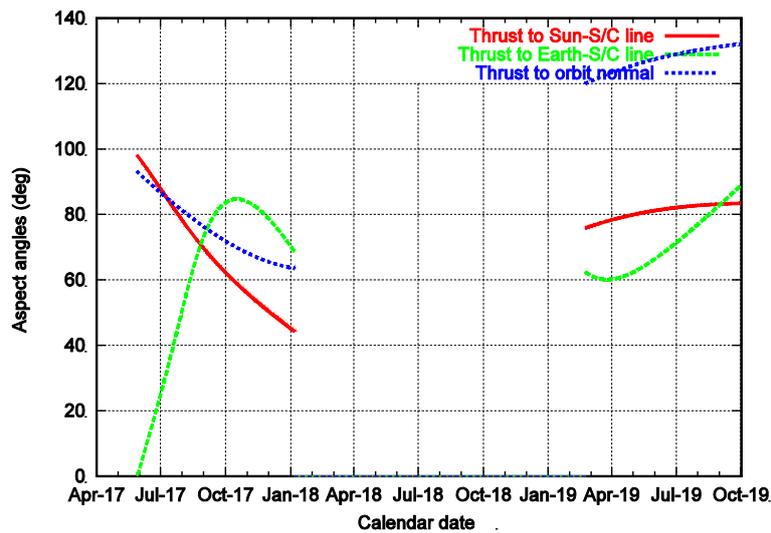
**Figure8**. Thrust aspect angles with respect to the Sun-S/C line, the Earth-S/C line and the orbit normal

A similar trajectory was found departing at mid June 2015 and with arrival to the asteroid in October 2017 but it was considered too challenging to meet the launch date.

## 3- IMPACT UNCERTAINTY COMPUTATION

The first task to be undertaken after rendezvous with the asteroid will be to perform an extensive tracking campaign aimed at shrinking the predicted uncertainty ellipsoid on the potential impact date in September 2022. On the expected rendezvous date at the end of September 2029 the asteroid will be favorably positioned with respect to the Earth and the Sun for tracking. Without performing any detailed analysis about orbit determination capability we assume a similar position and velocity accuracy to the one available for asteroid Itokawa at the time of rendezvous with the Hayabusa spacecraft at the end of 2005. Thanks to Arecibo radar measurements conducted in 2004 the asteroid velocity was known with an accuracy of about 0.2 mm/s (1-sigma). We use the same figures to construct a covariance matrix whose eigenvectors coincide with the unperturbed Frenet axes on November 1st 2019 and, in particular, are characterized by having the corresponding largest eigenvalues along the tangential direction as reasonable from an orbit determination point of view. In this way the covariance matrix expressed in unperturbed Frenet axes yields:

$$P = diag(\sigma_x^2, \sigma_y^2, \sigma_z^2, \sigma_{vx}^2, \sigma_{vy}^2, \sigma_{vz}^2),$$

where the standard deviations $\sigma_i$ of the position and velocity along the tangential, normal and binormal direction are given in Table 3.

**Table 3**. Expected 1-sigma post-rendezvous orbit determination errors for 2015$_{PDC}$

|  | position error | velocity error |
| --- | --- | --- |
| tangential (x) | 3.6 km | 0.2 mm/s |
| normal (y) | 150 m | 0.01 mm/s |
| binormal (z) | 400 m | 0.05 mm/s |

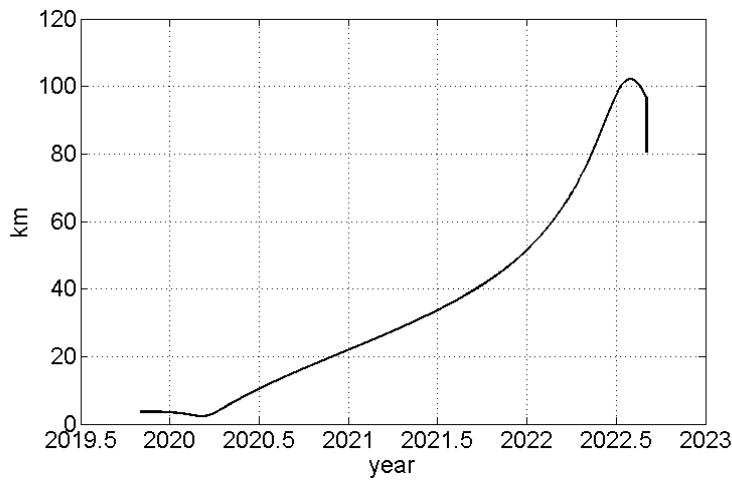

**Figure 9**. Propagated 1-sigma error ellipsoid largest semi-axis from rendezvous up to impact date

By linear propagation of the covariance matrix at the rendezvous date up to the expected impact date one would obtain an increase in the final position uncertainty up to about 80 km 1-sigma (Fig. 9) which is sufficient to confirm or rule out an impact. The projection of the covariance ellipsoid on the Earth surface would yield a 120x6 km size ellipse (1-sigma) with the longest axis roughly along the path of risk direction. That is enough to pinpoint the predicted impact location on the world map and decide about further action. In that regard, one should probably refer to a

600x30 km size 5-sigma uncertainty ellipse, which would cover a confidence region within $3.7 \times 10^{-6}$ probability of error.

It is important to add that, thanks to the continuous orbit determination campaign of the rendezvous spacecraft, the impact spot uncertainty will gradually diminish as the impact date approaches, finally reducing the potentially affected impact area to a circle corresponding to the radius of action of the impact (here assumed to be 100 km).

## 5- IMPACT POINT RETARGETING

Owing to the relatively short warning time and the complexity of the interplanetary trajectory to reach 2015 PDC it will be difficult to obtain, using non-nuclear method, a total deflection sufficient to have the asteroid missing the Earth by a safe amount. This is especially true if the asteroid diameter is confirmed to be in the range of 200-250 m or higher. On top of that, the risk of having the asteroid striking a different part of the planet as a result of a mission failure or an unexpected deflection outcome may in the end play against a full-scale deflection mission.

Consider for instance the case in which the asteroid is predicted to impact Northern India. A kinetic impact or slow push full deflection mission would cause the asteroid impact point to pass over at least five different countries according to the total momentum transmitted to the asteroid by the impulsive or slow push method. Slowing down the asteroid would cause its impact point to travel westward across India, Pakistan, Afghanistan, Iran and Turkey before an impact threat can be completely eliminated. On the other hand, by accelerating the asteroid path its impact point would travel eastward through Bangladesh, Myanmar, Thailand, Laos, Vietnam and Philippines.

Clearly, a full deflection mission should be designed with a large enough safety factor to make an unwanted impact totally unlikely. This may be very complicated, the more so the larger the asteroid.

An alternative mitigation action would be to shift the predicted impact point to a small amount just enough to have the asteroid striking a deserted or extremely low populated area to drastically reduce any expected damage or casualty risk. This is the approach to be analyzed in the following.

## 6- IMPACT DAMAGE MITIGATION ASSESSMENT

A land impact of a 150 to 500 m diameter asteroid may cause a significant amount of human casualties and damage to infrastructures, the magnitude of such effects depending on whether the impact point will be located close to a populated or urbanized area. Estimating such damage and its variation as a function of the impact point displacement is paramount in order to construct an effective mitigation action.

Assuming that the asteroid is accurately and continuously tracked by the rendezvous spacecraft during the whole deflection mission a rather accurate estimation, possibly within a few kilometers, of the impact point location will be available a few months before the impact. If the asteroid impact point is retargeted towards a minimum population density region there will still be the need to evacuate a sufficiently large area falling within the expected radius of action of the asteroid impact. The evacuation procedure will have a cost depending on the total number of inhabitants to be displaced. In addition, all infrastructures located within the radius of action of the impact may suffer irreversible damage and that will need to be quantified.

In light of the above considerations, two complementary mitigation indexes were considered. The first, referred as *Human Casualty Index (HCI)* is here defined as the ratio between the value (P) of the estimated population potentially affected by the impact in a generic location and the value ($P_0$) computed for the initial, undeflected impact location.

$$\text{HCI} = \frac{P}{P_0}$$

The total population P can be obtained from the 2015 projections of global population density data [7].

Note that the approach followed here differs from the one used by other authors (see for instance [6]) as no convolution with the impact position uncertainty distribution is applied. This is justified by the assumption that the impact uncertainty is much smaller than the critical radius of action of the asteroid impact damage.

The second mitigation index is called *Infrastructure Damage Index* (IDI) defined as the ratio of integrated nightime light intensity over the asteroid impact damage area computed for the deflected and undeflected impact point:

$$\text{HCI} = \frac{\int_A s dA}{\int_{A0} s dA}$$

Here, *s* is the light intensity per unit area quantized in 101 gray-scale levels in a high-resolution, nighttime image of the Earth provided by NASA [8].

The IDI index computed here is important as areas with similar population density can result in very different estimated damage figures depending on the presence of valuable infrastructures whose damage can never be completely eliminated.

## 5- DEFLECTION SCENARIOS

A simulation campaign to evaluate the IBS deflection capability of asteroid 2015 PDC has been conducted for each of the three worst case scenarios of Table 2 using the same high-fidelity orbit propagation model previously described. A continuous low-thrust deflection force is transmitted to the asteroid along or opposite to the instantaneous velocity vector. The thrust is applied starting in November 2019 (one month after rendezvous) and for a duration between one and 33 months obtaining a series of displaced impact points according to the duration of the deflection. As a reference, a 250 m diameter asteroid with 2 g/cm3 density has been perturbed with a 185 mN deflection force at 1 AU, with a varying magnitude inversely proportional to the distance from the Sun elevated to the 1.7 power.

The HCI and IDI are computed as a function of the shifted impact point during the whole deflection period. In this way it is possible to estimate the duration of deflection, which minimizes evacuation costs and infrastructural damage during the whole deflection period. The estimated population within the asteroid impact radius of action has also been computed.

**5.1 New Delhi impact**

A predicted impact of 2015 PDC near New Delhi would prompt to an immediate deflection action. By looking at Fig. 10-11 one immediately realizes that the nearest unpopulated region where the asteroid could be sent is central Afghanistan. Obviously, such a decision would be extremely delicate from a political point of view: the asteroid would need to travel across a highly populated region in Pakistan, who may strongly oppose such a move. However, political negotiations could finally lead to a solution in which the threatened country (India) decides to pay Pakistan a toll for the right of having the asteroid impact point passing over its territory and to pay a fee to Afghanistan for using its unpopulated territory to "absorb" the asteroid impact threat. India may decide to conduct the deflection mission on its own or have to pay another country or space agency for doing that. The deflecting country/agency may need to accept liability for any failure to accomplish its mission, which would contribute to the price to be paid in exchange for the deflection.

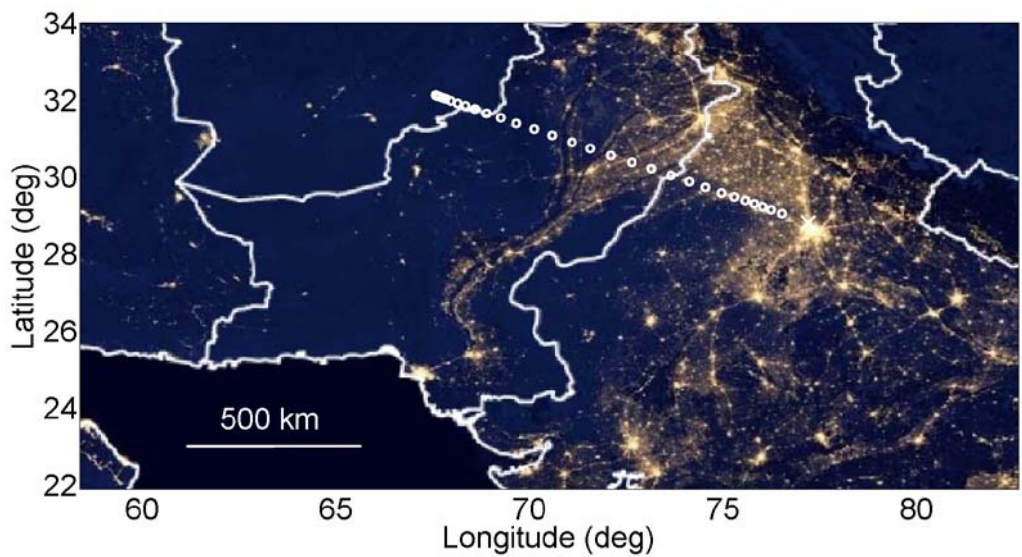

**Figure10**. IBS deflection track starting from a predicted impact in New Delhi. A cross marks the initial impact point.

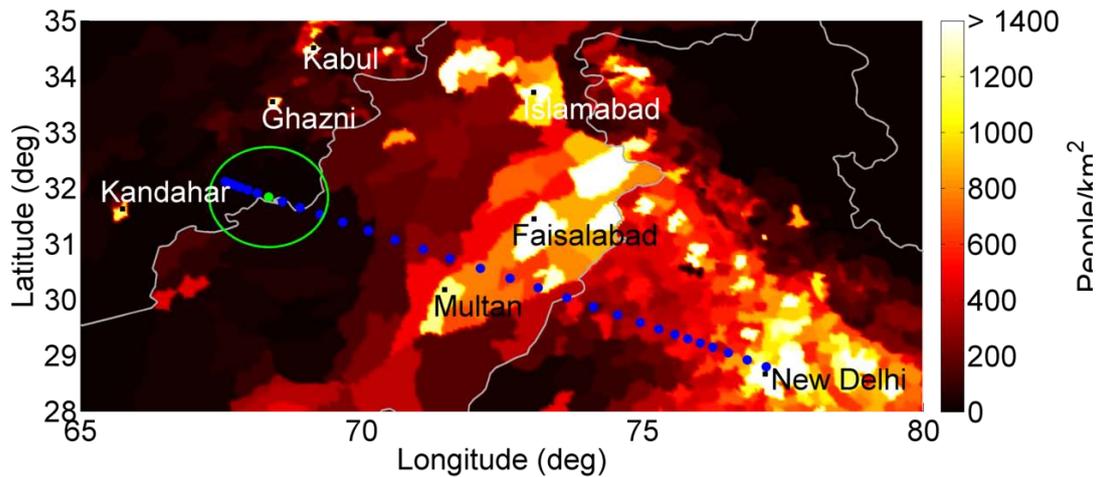

**Figure11**. Same as Fig 10 but starting from a population density map.

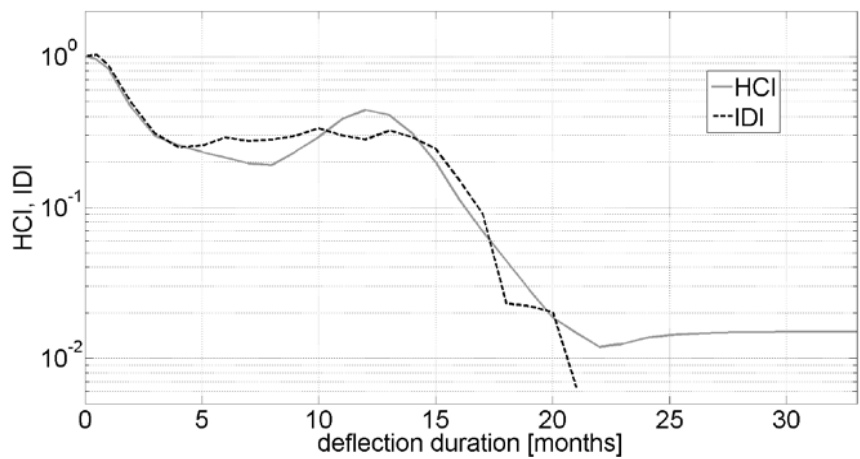

**Figure 12**. Variation of the HCI and IDI indexes with respect to the initial values during the deflection, starting from a predicted impact in New Delhi.

The shift of the impact point during the deflection is shown in Figure 10-11, in which the displaced impact points corresponding to each month of deflection are denoted with circles. The achieved shift is larger in the first half of the deflection period, and gets significantly smaller in the last months, when the applied deflection becomes less efficient due to the time before impact decreasing. Figure 12 highlights the benefit of such a deflection. It is possible to reduce the estimated casualties by two orders of magnitude and to virtually eliminate any infrastructural damage after 22 months of applied deflection (from that point on the integrated illumination map becomes zero). This corresponds to relocating the impact point in the sparsely populated Paktika province of Afghanistan, close to the border with Pakistan. There appear to be virtually no modern infrastructures in the area according to the nightlight intensity map, which means that the estimated damage would be limited to the cost of a coordinated evacuation of a relatively low number of inhabitants (below 500,000 according to Fig.19) and possible compensation for minor damages.

For shorter deflections (of less than 15 months) it is still possible to reduce casualties and infrastructural damage by roughly one order of magnitude and have the impact point moved to the Punjab region of Pakistan. This would have a reduced benefit in terms of infrastructural damage as reflected by a somewhat higher level of nighttime illumination compared to the secluded Paktika region.

## 5.2 Dhaka impact

A deflection action to avert an impact in Dhaka would be directed at displacing the impact point towards Myanmar, also crossing the Siliguri corridor belonging to India. Fig. 13-14 show the displaced impact points for each additional month of deflection duration. The optimal deflection, which reduces casualties and infrastructure damage by two orders of magnitude, is achieved after 13 months (see Fig. 15). This corresponds to shifting the impact point to the mountainous and scarcely inhabited Chin State in Myanmar. A considerably large number of people (about 1.2 million) would still fall in the estimated 100-km radius of action of the asteroid damage. However, one should also consider the benefit of having the asteroid landing in a mountainous region with higher damage attenuation capability.

As before, this would imply negotiations in which Bangladesh agrees to pay an impact fee to Myanmar for the right to displace the impact point to its territory, which would also cover the cost of evacuating the affected zones. Since some of the impact points would affect both Bangladesh and India, this could lead to an internationally coordinated action in which both states cooperate in the deflection.

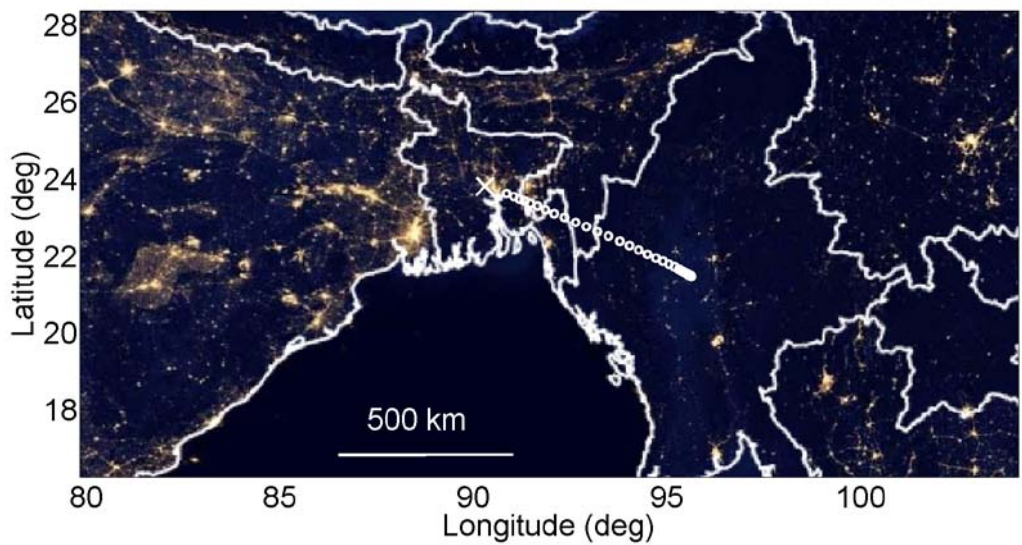

**Figure 13**. IBS deflection track starting from a predicted impact in Dhaka. A cross marks the initial impact point.

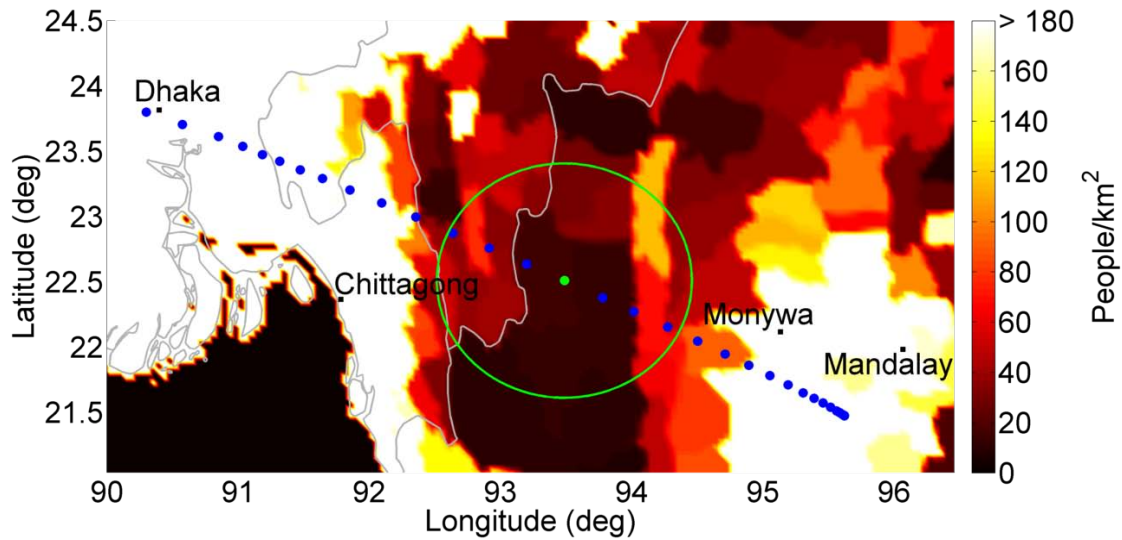

**Figure 14**. Same as Fig 13 but starting from a population density map.

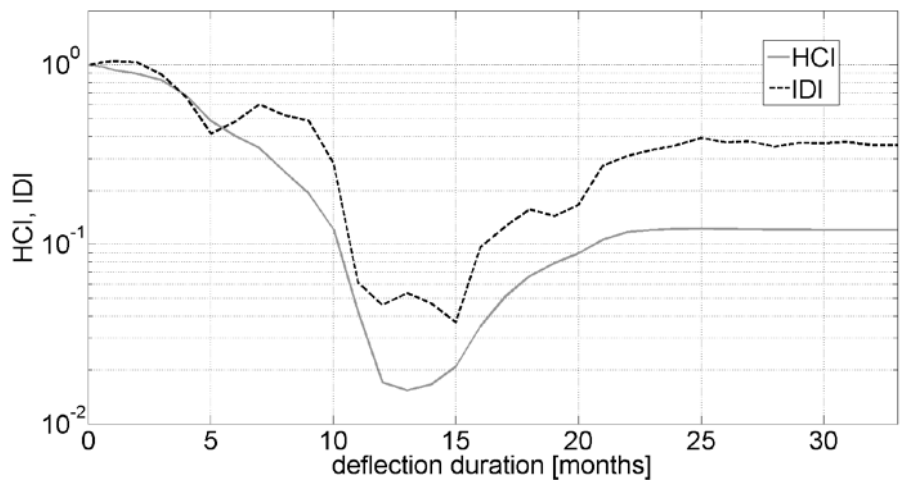

**Figure 15**. Variation of the HCI and IDI indexes with respect to the initial values during the deflection, starting from a predicted impact in Dhaka.

### 5.3 Tehran impact

A predicted impact in Tehran presents the easiest challenge, since it is possible to cut potential casualties by more than two orders of magnitude already after between one and two months of deflection. Fig. 16-17 reveal that the impact points which allow to minimize the damage are still contained in Iran, thanks to the presence of large swaths of deserted land that could accommodate the impact. This would greatly simplify political issues. However, the mission budget could still require a degree of international cooperation. In addition, international coordination would greatly benefit the transparency of the deflection action and provide opportunities for exchanging technical and scientific knowledge derived from the mission.

As shown in Fig. 18, after about 5 months of deflection both the HCI and IDI indexes start to rise again, with the IDI reaching a peak around month 8 due to the path of the deflected impact points passing through the city of Gonabad. This demonstrates the importance of accurately planning the deflection maneuver as to avoid causing collateral damage. Between month 12 and 16 the IDI becomes zero as a consequence of the asteroid landing in a completely dark zone, which would be ideal for the impact. The estimated evacuation action in that case would involve less than 170,000 inhabitants (Fig. 19).

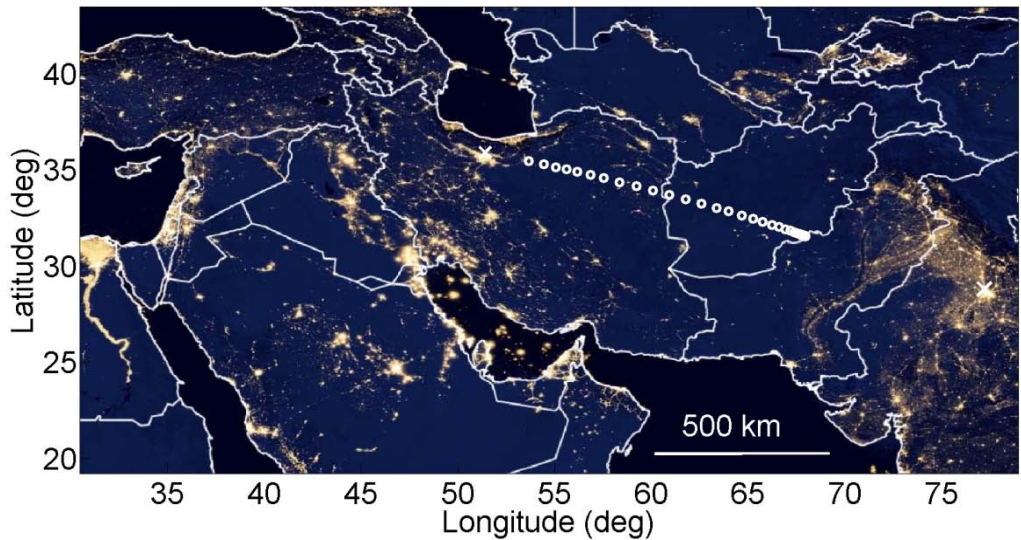

**Figure 16**. IBS deflection track starting from a predicted impact in Teheran.

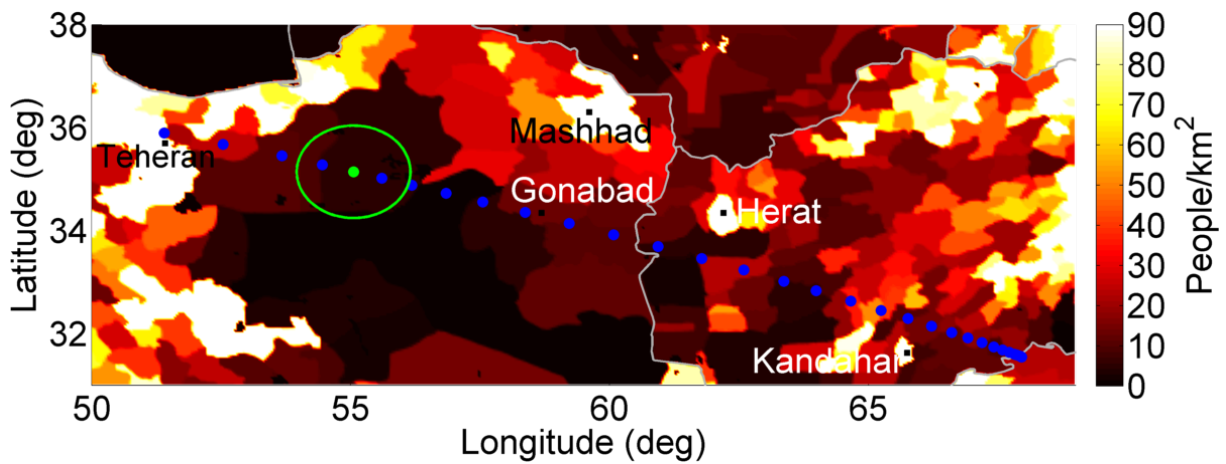

**Figure 17**. Same as Fig 16 but starting from a population density map.

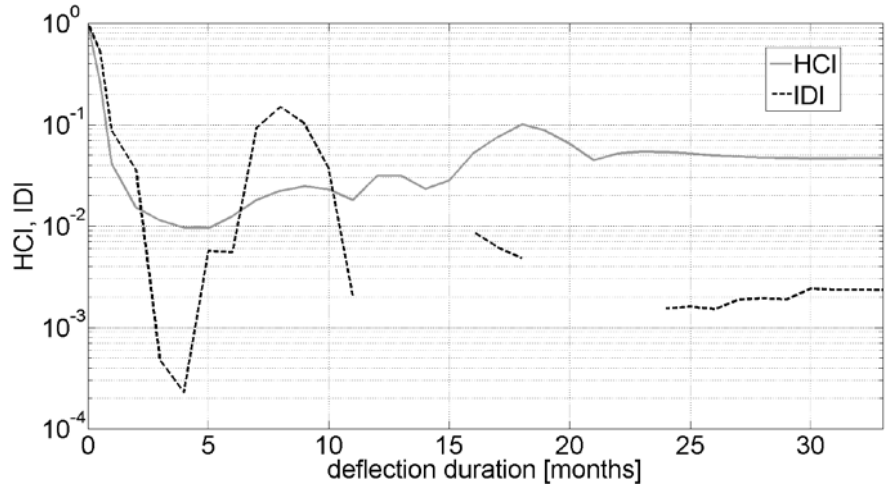

**Figure 18**. Variation of the HCI and IDI with respect to the initial values during the deflection, starting from a predicted impact in Tehran.

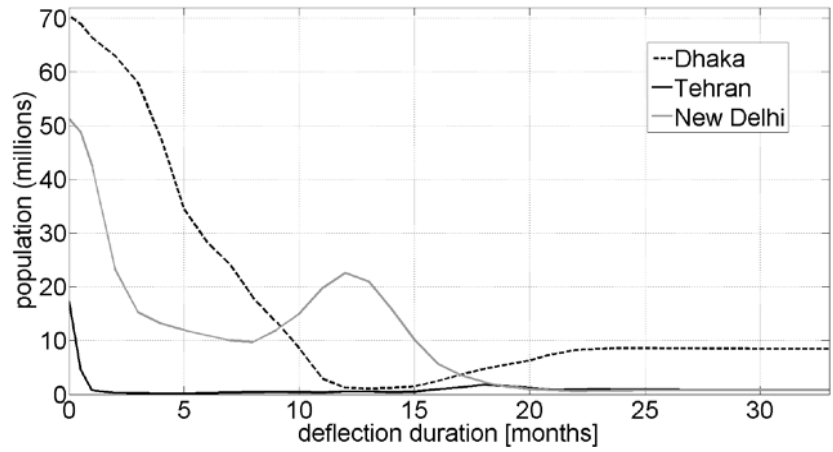

**Figure 19**. Total population to be evacuated in case of an impact as a function of the deflection duration.

## 6- VERY LARGE ASTEROID IMPACT

The presented mission design concept was based on a maximum assumed asteroid diameter of 250 m, allowing a relatively cheap rendezvous and deflection mission with an estimated 1200 kg spacecraft mass at interplanetary injection. An asteroid of 500 m in diameter, corresponding to the maximum size considered by the conference organizers, would require a deflection capability 8 times as large. A single IBS spacecraft with such capability exceeds current technology limitations and would probably not be selected as baseline design for the first rendezvous mission. Multiple IBS could be flown later and act in parallel to provide the required deflection capability. Collision avoidance between multiple deflection units would certainly complicate the mission. Another option is to send a number of kinetic impactor spacecraft following the IBS rendezvous, providing a series of velocity changes whose effect would be measured by the IBS spacecraft. The latter would provide the last deflection refinement to have the asteroid falling in a minimum damage area as discussed in the previous sections. An investigation of the effectiveness of a kinetic impacting deflection method for 2015 PDC has been presented in this conference [9]. Finally a single deflection using a nuclear explosion could be considered as an alternative solution having the advantage of providing full deflection of even the largest expected asteroid. The risk and cost of such mission would need to be carefully examined.

## CONCLUSIONS

A preliminary design and analysis of a slow push deflection mission of fictitious asteroid 2015 PDC has been conducted making use of an ion beam shepherd(IBS) spacecraft. The IBS system could fit, with reduced impact in mass and complexity, in a small 1200 kg spacecraft, which would be sent to rendezvous with the asteroid and perform an indispensable refinement of its orbit accuracy to confirm or rule out an impact. In case an impact near a major city is predicted it is proposed to employ the IBS to accurately displace the impact point towards a deserted or very low populated neighboring region, which greatly reduces mission costs when compared to full deflection. Two indexes are introduced to quantify the benefit of such action by looking at potential population and infrastructure damage. The approach appears to be feasible also from the political point of view as long as mission costs, liabilities and compensations are agreed upon at international level.

## ACKNOWLEDGEMENTS

The research leading to these results has received funding from the European Union Seventh Framework Programme (FP7/2007-2013) under grant agreement 317185 (Stardust) and by the Spanish Ministry of Economy and Competitiveness within the framework of the research project "Dynamical Analysis, Advanced Orbital Propagation, and Simulation of Complex Space Systems" (ESP2013-41634-P).
Last but not least, we would like to thank the organizers of the 2015 Planetary Defense Conference for proposing this very interesting problem.